\input amstex
\documentstyle{amsppt}
\magnification 1200
\NoBlackBoxes
\vcorrection{-9mm}
\input epsf

\topmatter
\title
           On osculating framing of real algebraic links
\endtitle
\author
          Grigory Mikhalkin and Stepan Orevkov
\endauthor
\abstract
For a real algebraic link in $\Bbb{RP}^3$,
we prove that its encomplexed writhe (an invariant introduced by Viro) is maximal
for a given degree and genus if and only if its self-linking number with respect to the
framing by the osculating planes is maximal for a given degree.
\endabstract
\thanks
Research is supported in part by  the SNSF-grants
159240, 159581, 182111  and the NCCR SwissMAP project  (G.M.),
and by RFBR grant No 17-01-00592a (S.O.)
\endthanks

\address        Universit\'e de Gen\`eve, Section de Math\'ematiques,
                Battelle Villa, 1227 Carouge, Suisse.
\endaddress
\email          grigory.mikhalkin\@unige.ch
\endemail

\address
                Steklov Mathematical Institute, Gubkina 8, 119991, Moscow, Russia;
\smallskip
                IMT, Universit\'e Paul Sabatier,
                118 route de Narbonne, 31062, Toulouse, France.
\endaddress
\email          orevkov\@math.ups-tlse.fr
\endemail
\endtopmatter

\def\R{\Bbb R}

\def\P{\Bbb P}    \let\pp=\P
\def\RP{\Bbb{RP}} \let\rp=\RP
\def\CP{\Bbb{CP}}
\def\osc{\operatorname{osc}}
\def\lk{\operatorname{lk}}
\def\ps{\operatorname{ps}}
\def\hdp{h}                 

\def\refKlein     {1}
\def\refMOdan     {2}
\def\refMOmw      {3}
\def\refM         {4}
\def\refSchuh     {5}
\def\refViroPluck {6}
\def\refViroEW    {7}

\def\eqBB       {1}
\def\MainIneq   {2}
\def\eqPos      {3}
\def\eqLK       {4}
\def\eqOscK     {5}
\def\ineqLK     {6}
\def\ineqG      {7}

\def\figTL   {1}

\def\lemFlex  {2.1}
\def\remFlex  {2.2}

\document

\head   1. Introduction and statement of the main result
\endhead

By {\it real algebraic curve} in $\RP^3$ we mean
a complex curve in $\CP^3$ invariant under complex conjugation. We use the
same notation for a real curve and the set of its complex points and,
if it is denoted by $A$, then $\R A$ stands for the set of real points
which is called a {\it real algebraic link} if it is non-empty and $A$ is smooth.
A real algebraic link is called {\it maximally writhed} or {\it $MW_\lambda$-link}
if $|w_\lambda(L)|$ (a variation of Viro's invariant [\refViroEW])
attains the maximal possible value $(d-1)(d-2)/2-g$ where $d$ and $g$
is the degree and genus of $A$ respectively. 
We refer to [\refMOmw] for a precise definition of $w_\lambda$.

In [\refMOmw, Thm.~2] we proved that several topological and
geometric invariants are maximized on $MW_\lambda$-links.
In this paper we add one more item to this collection: we show that
the self-linking number of $L$ with respect to the osculating framing
attains its maximal value (for links of a given degree) if and only if
$L$ is an $MW_\lambda$-link. The proof is very similar to that of the main theorem of
[\refMOmw]. Let us give precise definitions and statements.

Let $L$ be an oriented link in a rational homology 3-sphere. A {\it framing} of $L$ is a
continuous 1-dimensional subbundle of the normal bundle of $L$ or, equivalently,
a continuous field (defined on $L$) of 2-dimensional planes tangent to $L$.
Given a framed oriented link $L$, its {\it self-linking number} is defined as follows.
Let $F$ be an embedded annulus or M\"obius band with core $L$, tangent to the
framing. Then the self-linking number is $\frac12\lk(L,\partial F)$ where
the boundary $\partial F$ of $F$ is oriented so that $[\partial F]=2[L]$ in $H_1(F)$.

For an oriented link $L$ in $\RP^3$, the {\it osculating framing} is the framing defined by
the field of osculating planes.
We denote the self-linking number of $L$ with respect to this
framing by $\osc(L)$. If $L$ is a non-oriented link and $O$ an orientation of $L$, we
use the notation $\osc(L,O)$ which is self-explained.

Recall that a smooth irreducible real algebraic curve $A$ is called an {\it $M$-curve} if $\R A$
has $g+1$ connected components where $g$ is the genus of $A$. In this case
$\R A$ divides $A$ into two halves. The boundary orientation on $\R A$ induced
by any of these halves is called a {\it complex orientation}.
The main result of the paper is the following.

\proclaim{ Theorem 1 }
Let $L=\R A$ be an irreducible real algebraic link of degree $d\ge 3$
and $O$ be an orientation of $L$.
Then:
\smallskip
\roster
\item"(a)"
          $|\osc(L,O)|\le d(d-2)/2$.
\smallskip
\item"(b)"
          $|\osc(L,O)|=d(d-2)/2$ if and only if $L$ is an $MW_\lambda$-link
          {\rm (by [\refMOmw, Thm.~2], in this case $A$
          is an $M$-curve of genus at most $d-3$)} and
          $O$ is its complex orientation.
\endroster
\endproclaim

\noindent{\bf Remark. } In the space of real algebraic links of a given degree and genus
we can distinguish three kinds of ``walls". The walls of the first kind correspond to
curves with a double point with real local branches. When crossing such walls,
both invariants $w_\lambda(L)$ and $\osc(L)$ are changed by $\pm2$.
The walls of the second kind correspond to curves with a real double point with
complex conjugate local branches. When crossing such walls, $w_\lambda(L)$ does change
but $\osc(L)$ does not. The third kind of wall corresponds to curves
which have a local branch parametrized by $t\mapsto(t,t^3+o(t^3),t^4+o(t^4))$ in
some affine chart. When crossing such a wall, $w_\lambda(L)$ does not change
but $\osc(L)$ does. So, in general, the invariants $w_\lambda(L)$ and $\osc(L)$
are more or less independent. Nevertheless, Theorem 1 implies that
the chamber where they have maximal value is bounded only by the walls of
the first kind -- common for the both invariants.


\head     2. A variant of Klein's formula for the number of real inflection points
\endhead

Let $C\in\P^2$ be a nodal real irreducible algebraic curve.
It may have three types of nodes: real nodes with real local
branches of $C$, real nodes with imaginary local branches
of $C$, or non-real nodes (coming in conjugate pairs).
Denote the number of nodes of each type with $h$, $e$, and $i$
respectively.

A {\it real flex} is a local real branch of $C$ with the
order of tangency $\omega$ to its tangent line greater than 1
(i.e. the local intersection number is $\omega+1\ge 3$).
The multiplicity of a real flex is $\omega-1$. In an affine
chart of $\pp^2$ a flex corresponds to a critical point
of the Gauss map. It is easy to see that the multiplicity
of a flex equals to the multiplicity of the corresponding
critical point. Thus a multiple flex can be thought of as
$\omega-1$ ordinary flexes collected at the same point.
We denote with $F$ the number of flexes counted with
multiplicities.

A {\it solitary real bitangent} is a real line $L\subset\pp^2$
which is tangent to $C$ at a non-real point (and thus
also at the complex conjugate point). The {\it multiplicity}
of $L$ is the sum of the orders $\omega$ over all local
branches of $C\setminus\RP^2$ tangent to $L$.
We denote with $B$ the number of solitary real bitangents
counted with multiplicities. Clearly, $B$ is an even number.

\proclaim{Lemma \lemFlex} {\rm(Klein's formula [\refKlein] for nodal curves).}
For a nodal real irreducible curve of degree $d$ in $\pp^2$ we have
$$
      F+B=d(d-2)-2h-2i.
$$
\endproclaim

\demo{ Proof }
As in [\refViroPluck], we use additivity of the Euler characteristic $\chi$
to derive Klein's formula. Let $\nu:\tilde C\to C$ be the normalization.
The space of all real lines in $\pp^2$ is homeomorphic to $\RP^2$, and
thus has the Euler characteristic 1. For a real line $L$ the set
$\nu^{-1}(L)$ consists of $d$
distinct points unless $L$ is tangent to $C$. Each tangency decreases
the size of this set by $\omega$.

Consider the space $X=\{(p,L)\mid p\in C,\ L\ni p\}$,
where $L\subset\rp^2$ is a real line.
From the observation above we deduce
$$
       \chi(X)+B+F+\chi(\R \tilde C)=d.
$$
Note that $\chi(\R \tilde C)=0$ and
$\chi(X)=\chi(\nu^{-1}(C\setminus\R C))=\chi(\tilde C)-2e$,
as each point of $\R C$ lifts to a circle in $X$ while $\chi(S^1)=0$.
The lemma now follows from the adjunction formula
$\chi(\tilde C)=3d-d^2+2e+2h+2i$.
\enddemo

\noindent
{\bf Remark \remFlex.}
Lemma \lemFlex\ can be also obtained as an almost immediate consequence
from Schuh's generalization [\refSchuh] of another Klein's formula
$$
    d - \sum_{x\in C\cap \RP^2}(m(x)-r(x)) =
    d^\vee - \sum_{x\in C^\vee\cap {\RP^2}^\vee}(m^\vee(x)-r^\vee(x))
$$
(see [\refViroPluck, Thm. 6.D] for a proof via Euler characteristics) combined with the
class formula $d^\vee=d(d-1)-2e-2h-2i$. Here $C^\vee$ is the dual curve, $d^\vee$ is its degree, $m(x)$ and $r(x)$ (resp. $m^\vee(x)$ and $r^\vee(x)$)
are the multiplicity and the number of real local branches of $C$ (resp. of $C^\vee$)
at $x$.


\head     3. Proof of the main theorem
\endhead

Let $L=\R A$ be a smooth irreducible real algebraic link
of degree $d$ endowed with an orientation $O$.
Let $\Cal U$ be the set of points $p$ in $\RP^3\setminus L$ such that
the projection of $L$ from $p$ is a nodal curve.

Fix a point $p\in \Cal U$.
Let $C_p=\pi_p(A)$ where $\pi_p:\P^3\setminus\{p\}\to\P^2$ is the linear projection 
from $p$.
Consider the field of tangent planes to $L$ passing through $p$,
(so-called blackboard framing).
Let $b_p(L)$ be the self-linking number with respect to it.
We have
$$
     b_p(L) = \sum_q s(q), \qquad\text{thus}\qquad |b_p(L)| \le \hdp(C_p)    \eqno(\eqBB)
$$
where $q$ runs the hyperbolic (i.~e., with real local branches)
double points of $C_p$, $\hdp(C_p)$ is the
number of them, and $s(q)$ is the sign
of the crossing at $q$ in the sense of knot diagrams.
The difference $|\osc(L)-b_p(L)|$ is bounded by one half of the number of those points
where the osculating plane passes through $p$. This is the
number of real flexes of $C_p$ which we denote by $f(C_p)$.
We have $f(C_p) \le d(d-2) - 2\hdp(C_p)$ by Lemma \lemFlex.
Thus
$$
     |\osc(L)| 
          \le |\osc(L) - b_p(L)| + |b_p(L)|
          \le\tfrac12 f(C_p) + \hdp(C_p) \le \tfrac12 d(d-2)          \eqno(\MainIneq)
$$
which is Part (a) of Theorem 1.

Now suppose that $|\osc(L)|=d(d-2)/2$. Then for any choice of $p\in\Cal U$
we have the equality sign every\-where in (\MainIneq),
in particular, we have the equality sign in (\eqBB), i.e.,
all crossings are of the same sign, say, positive:
$$
        s(q)=+1 \qquad\text{for any hyperbolic crossing $q$ of $C_p$.}  \eqno(\eqPos)
$$

By Lemma \lemFlex, the equality sign in the last inequality of (\MainIneq)
implies that all flexes of $C_p$ are ordinary for any choice of $p\in\Cal U$.
This implies that $L$ has non-zero torsion at each point. Indeed, otherwise
there exists a real plane $P$ which has tangency with $L$ of order greater than $3$.
It is easy to check that $\Cal U$ has non-empty intersection with any plane, thus
we can choose a point $p\in\Cal U\cap P$, and then $C_p$ would have a $k$-flex with $k>3$.
Moreover, the positivity of all crossings for any generic projection implies that the
torsion is everywhere positive (cf. the proof of [\refMOdan, Prop.~1]).

Similarly to [\refMOdan, \refMOmw], we derive from these conditions that the
real tangent surface $TL$ (the union of all real lines in $\RP^3$ tangent to $L$)
is a union of (non-smooth) embedded tori.
Indeed, suppose that two tangent lines cross.
Let $P$ be the plane passing through them (any plane passing through them if they coincide)
and let $\ell$ be the line passing
through the two tangency points.
Let $p$ be a generic real point on $\ell$. Then $C_p$ has two real local
branches at the same point such that each of them is either singular or
tangent to the line $\pi_p(P)$. Since $L$ has non-zero torsion, all singular
branches of $C_p$ are ordinary cusps.
Then we can find a generic point close to $p$ such that the projection from it
does not satisfy (\eqPos).

Let $K_1,\dots,K_n$ be the connected components of $L$, and let $TK_i$
be the connected component of $TL$ that contains $K_i$
(the union of real lines tangent to $K_i$).
The same arguments as in [\refMOmw, Lemma 4.12] show that, for some positive
integers $a_1,\dots,a_n$, there exist real lines
$\ell_i$, $\ell'_i$, $i=1,\dots,n$, such that (for suitable choice of
the orientations) the linking numbers of
their real loci $l_i=\R\ell_i$ and $l'_i=\R\ell'_i$ with the components
of $L$ are:
$$
     2\lk(l_i,K_i)=a_i+2, \qquad 2\lk(\l'_i,K_i)=a_i.             \eqno(\eqLK)
$$
Moreover, each $TK_i$ splits $\RP^3$ into two solid tori $U_i$ and $V_i$
such that $l_i\subset U_i$, $l'_i\subset V_i$, the homology classes
$[l_i]_U$ and $[l'_i]_V$ generate $H_1(U_i)$ and $H_1(V_i)$ respectively, and
we have $[K_i]_U=a_i[l_i]_U$ and $[K_i]_V=(a_i+2)[l'_i]_V$.
It follows that
$$
    2\osc(K_i)=a_i(a_i+2)                                      \eqno(\eqOscK)
$$
(the linking number of $K_i$ with its small shift disjoint from $TL$).
Indeed, if $K_i$ is parametrized by $t\mapsto r(t)$ and the torsion is
non-zero, then $TK_i$ has a cuspidal edge along $K_i$ and a small shift of $K_i$
in the direction of the vector field
$\ddot r$ is disjoint from $TK_i$ (see Figure \figTL). A priori this argument
proves (\eqOscK) up to sign only. However the positivity of the torsion implies
that $\osc(K_i)$ is positive.

\midinsert
\centerline{\epsfxsize=20mm \epsfbox{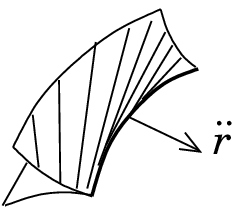}}
\botcaption{Figure \figTL}
\endcaption
\endinsert

If $L$ is connected (i.~e., $n=1$), it remains to note that
then the condition $2\osc(K_1)=d(d-2)$ implies $(a_1+2)a_1=d(d-2)$, hence
$a_1=d-2$. Thus $L$ satisfies
Condition (v) of [\refMOmw, Thm.~1] which concludes
the proof that $L$ is an $MW$-knot.

If $L$ is not necessarily connected, we argue as follows.
By Murasugi's result [\refM, Prop.~7.5] (see also [\refMOmw, Prop.~1.2]),
the number of crossings of any projection of $K_i$ is at least
$(a_i+2)(a_i-1)/2$. Hence, for $\hdp=\hdp(C_p)$, we have
$$
   2\hdp \ge \sum_{i=1}^n(a_i+2)(a_i-1) + \sum_{i\ne j}|\lk(K_i,K_j)|.  \eqno(\ineqLK)
$$
On the other hand, if we choose $p$ on a line passing through a pair of complex
conjugate points, then $C_p$ has at least one elliptic double point (i.~e.,
a real double point with complex conjugate local branches), whence by the
genus formula we obtain
$$
   \hdp \le (d-1)(d-2)/2 - g - 1 \le (d-1)(d-2)/2 - n                 \eqno(\ineqG)
$$
(the second inequality in (\ineqG) is the Harnack's bound). Hence
$$
\xalignat2
   d(d-2)&=2\osc(L) = 2\sum_{i=1}^n \osc(K_i) + \sum_{i\ne j}\lk(K_i,K_j) \\
       &\le \sum_{i=1}^n a_i(a_i+2) + 2\hdp - \sum_{i=1}^n(a_i+2)(a_i-1)
                      &&\text{by (\eqOscK) and (\ineqLK)} \\
       &=  2\hdp+2n+\sum_{i=1}^n a_i
        \le (d-1)(d-2)  +\sum_{i=1}^n a_i. &&\text{by (\ineqG)}
\endxalignat
$$
Thus $\sum a_i \ge d-2$ and we conclude that $L$ is an $MW_\lambda$-link.
This fact follows from
[\refMOmw, Prop.~1.1] (which implies that $\ps(L)=\sum a_i$)
combined with [\refMOmw, Thm.~2] (which claims, in particular,
that $L$ is an $MW_\lambda$-link as soon as $\ps(L)\ge d-2$).
Here we denote with $\ps(L)$ the {\it plane section number} of $L$. It
is a topological invariant of a link in $\RP^3$ defined in [\refMOmw] as
the minimal number of intersection points with a generic plane where the
minimum is taken over the isotopy class of the link.

Let us show that $O$ is a complex orientation of $L$. It is easy to see that
the plane section number is at most $d-2$ for any algebraic link of degree $d$.
Indeed, it is enough to consider a small shift of a non-osculating tangent plane
in a suitable direction. Thus the inequality in $\ps(L)=\sum a_i\ge d-2$
is in fact an equality. It follows that the equality is attained in all the
inequalities used in the proof, in particular, we have $|\lk(K_i,K_j)|=\lk(K_i,K_j)$
for $i\ne j$. Since all components of an $MW_\lambda$-link endowed with a
complex orientation are positively linked (see [\refMOmw]), we are done.
This completes the proof of the ``only if\," part of (b).

To prove the ``if\," part of (b), we notice that by
[\refMOmw, Thm. 3 and \S4.4], any $MW_\lambda$-link $L$
of degree $d$ and genus $g$ is a union of $g+1$ knots $K_0\cup\dots\cup K_g$
and  $\lk(K_i,K_j)=a_i a_j$, $i\ne j$, for some positive integers $a_0,\dots,a_g$
with $a_0+\dots+a_g=d-2$. Furthermore, the torsion of $L$ is everywhere positive
and each knot $K_i$ is arranged on its tangent surface $TK_i$ as described above,
thus (\eqOscK) holds for each $i$, and we obtain
$$
\split
   2\osc(L) &= \sum_{i=0}^g\osc(K_i)+\sum_{i\ne j}\lk(K_i,K_j)
            = \sum_{i=0}^g a_i(a_i+2) + \sum_{i\ne j} a_i a_j
\\&
            = \Big(\sum a_i\Big)^2 + 2\sum a_i = (d-2)^2+2(d-2) = d(d-2).
\endsplit
$$


\enddocument